\documentclass[11pt]{article}
\pdfoutput=1
\usepackage[T1]{fontenc}
\usepackage{amsfonts,amssymb,amsmath,amsthm,mathabx}

\usepackage{fullpage,xcolor}
\usepackage{mathabx}
\usepackage{graphicx}
\graphicspath{{./},{./images/},}

\usepackage{pdfpages,hyperref}
\hypersetup{colorlinks=true, linkcolor=blue,  anchorcolor=blue, citecolor=blue, filecolor=blue, menucolor=blue, pagecolor=blue, urlcolor=blue}

\usepackage{fancyhdr}
\newcommand{\myitemsymb}{$\blacktriangleright$}%\triangleright$}
\begin{document}

\thispagestyle{empty}

\ 
\vspace{4cm}
\begin{center}
{\Huge \bf
Proceedings of the third\\
``international Traveling Workshop\\
on Interactions between Sparse models\\ and Technology''\\[2mm]
\textsc{iTWIST'16}\\[.5cm]
}
{\Large Aalborg, Denmark}\\[2mm]
{August 24-26, 2016.}\\[1cm]
\includegraphics[width=.7\textwidth]{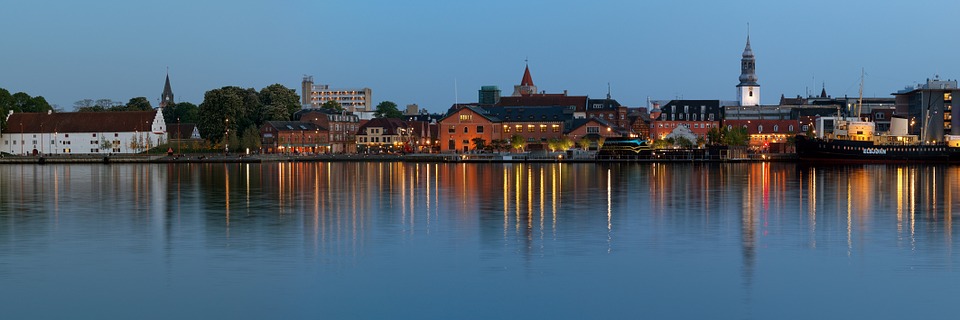}  
\end{center}

%%%%%%%%%%%%%%%%%%%%%%%%%%%%%%%%%%%%%%%%%%%%%%%%%%%%%%%%%%%%%%%%%%%%%%%%%%%%%%%%%%%%%%%%%%%%%%%%%%%
\newpage
\setcounter{page}{1}
\renewcommand{\thepage}{\roman{page}}
\ \\[1cm]
\begin{center}
\LARGE i\textsc{TWIST'16 Presentation}  
\end{center}

The third edition of the
``\emph{international - Traveling Workshop on Interactions between Sparse
models and Technology}'' (iTWIST) will take place in Aalborg, the 4th
largest city in Denmark situated beautifully in the northern part of
the country. The workshop venue will be at the Aalborg University
campus.

One implicit objective of this biennial workshop is to foster collaboration
between international scientific teams by disseminating ideas through
both specific oral/poster presentations and free discussions.
 
For this third edition, iTWIST’16 gathers about 50 international participants and features 8 invited talks,
12 oral presentations, and 12 posters on the following themes, all related to the theory, application and
generalization of the ``sparsity paradigm'':
\begin{itemize}
\setlength\itemsep{.1em}
\item[\myitemsymb] Sparsity-driven data sensing and processing (e.g., optics, computer vision, genomics, biomedical, digital communication, channel estimation, astronomy)
\item[\myitemsymb] Application of sparse models in non-convex/non-linear inverse problems (e.g., phase retrieval, blind deconvolution, self calibration)
\item[\myitemsymb] Approximate probabilistic inference for sparse problems
\item[\myitemsymb] Sparse machine learning and inference
\item[\myitemsymb] ``Blind'' inverse problems and dictionary learning
\item[\myitemsymb] Optimization for sparse modelling
\item[\myitemsymb] Information theory, geometry and randomness
\item[\myitemsymb] Sparsity? What's next?
  \begin{itemize}
  \item[\myitemsymb] Discrete-valued signals
  \item[\myitemsymb] Union of low-dimensional spaces,
  \item[\myitemsymb] Cosparsity, mixed/group norm, model-based, low-complexity models, ...
  \end{itemize}
\item[\myitemsymb] Matrix/manifold sensing/processing (graph, low-rank approximation, ...)
\item[\myitemsymb] Complexity/accuracy tradeoffs in numerical methods/optimization
\item[\myitemsymb] Electronic/optical compressive sensors (hardware)
\end{itemize}

\newpage
\ \\[1cm]
\begin{center}
\Large \textsc{Scientific Organizing Committee}  
\end{center}
\begin{itemize}
\item[\myitemsymb] {\bf Thomas Arildsen} (general chair)\\
  Department of Electronic Systems at Aalborg University, Denmark.
\item[\myitemsymb] {\bf     Morten Nielsen} (co-chair)\\
  Department of Mathematical Sciences at Aalborg University, Denmark.
\item[\myitemsymb]  {\bf    Laurent Jacques}\\
  Institute of Information and Communication Technologies,\\
 \quad Electronics and Applied Mathematics (ICTEAM) at Universit\'e
 catholique de Louvain, Belgium.
\item[\myitemsymb]  {\bf    Sandrine Anthoine}\\
  Marseille Institute of Mathematics, at Aix-Marseille Universit\'e
  and CNRS, France.
\item[\myitemsymb]   {\bf   Yannick Boursier}\\
  Center for Particle Physics of Marseille, at Aix-Marseille
  Universit\'e and CNRS, France.
\item[\myitemsymb]   {\bf   Aleksandra Pizurica}\\
  Department of Telecommunications and Information Processing at Ghent
  University, Belgium.
\item[\myitemsymb]   {\bf   Pascal Frossard}\\
  Electrical Engineering Institute at Ecole Polytechnique F\'ed\'erale
  de Lausanne, Switzerland.
\item[\myitemsymb]  {\bf    Pierre Vandergheynst}\\
  Electrical Engineering Institute at Ecole Polytechnique F\'ed\'erale de Lausanne, Switzerland.
\item[\myitemsymb] {\bf     Christine De Mol}\\
  Department of Mathematics and ECARES (European Center for Advanced
  Research in Economics and Statistics) at Brussels Free University, Belgium.
\item[\myitemsymb]  {\bf    Christophe De Vleeschouwer}\\
  Institute of Information and Communication Technologies,\\
 Electronics
  and Applied Mathematics (ICTEAM) at Universit\'e catholique de Louvain, Belgium.
\end{itemize}
  
% \noindent\textbf{Sponsors:} The iTWIST'14 organizing committee thanks the following sponsors for
% their help and fundings.
% \begin{center}
% \includegraphics[width=\textwidth]{sponsors_itwist14}
% \end{center}
\newpage

\section*{Table of contents}

\begin{itemize}
\footnotesize
%%%%%%%%%%%%%%
%main_02.tex:
\item
  \hyperlink{main_02.1}{
{\large ``Sparse matrix factorization for PSFs field estimation''}}\dotfill{\scriptsize (p.~\pageref*{pdf:main_02})}
\vspace{2mm} \newline F. Ngol\`e (CEA Saclay, France) and J.-L. Starck (CEA Saclay, France).
%--
%main_03.tex
\item \hyperlink{main_03.1}{\large 
``Sparse BSS with corrupted data in transformed domains''.}\dotfill{\scriptsize (p.~\pageref*{pdf:main_03})}
\vspace{2mm} \newline
C. Chenot (IRFU, CEA, France) and J. Bobin (IRFU, CEA, France). 
%--
%main_04.tex 
\item \hyperlink{main_04.1}{\large
``Randomness is sufficient! Approximate recovery from CS samples''.}\dotfill{\scriptsize (p.~\pageref*{pdf:main_04})}
\vspace{2mm} \newline
V.~Abrol (SCEE, IIT Mandi, India), P.~Sharma (SCEE, IIT Mandi, India) and Anil~K. Sao (SCEE, IIT Mandi, India).
%--
%main_05.tex
\item \hyperlink{main_05.1}{\large
``The Effect of Atom Replacement Strategies on Dictionary Learning''.}\dotfill{\scriptsize (p.~\pageref*{pdf:main_05})}
\vspace{2mm} \newline
Paul Irofti (U. Politehnica Bucharest, Romania).
%--
%main_06.tex
\item \hyperlink{main_06.1}{\large
``Blind Deconvolution of PET Images using Anatomical Priors''.}\dotfill{\scriptsize (p.~\pageref*{pdf:main_06})}
\vspace{2mm} \newline
St\'ephanie Gu\'erit (UCLouvain, Belgium),
    Adriana Gonz\'alez (UCLouvain, Belgium)
    (UCLouvain, Belgium), Anne Bol (UCLouvain, Belgium), John A. Lee (UCLouvain, Belgium) and Laurent Jacques (UCLouvain, Belgium).
%--
%main_07.tex
\item \hyperlink{main_07.1}{
\large
``A Non-Convex Approach to Blind Calibration for Linear Random Sensing Models''.}\dotfill{\scriptsize (p.~\pageref*{pdf:main_07})}
\vspace{2mm} \newline
Valerio Cambareri (UCLouvain, Belgium) and
    Laurent Jacques (UCLouvain, Belgium).
%--
%main_08.tex
\item \hyperlink{main_08.1}{\large
``Sparse Support Recovery with $\ell_\infty$ Data Fidelity''.}\dotfill{\scriptsize (p.~\pageref*{pdf:main_08})}
\vspace{2mm} \newline
K\'{e}vin Degraux (UCLouvain, Belgium),
    Gabriel Peyr\'{e} (CNRS, Ceremade, U. Paris-Dauphine, France),
    Jalal M. Fadili (ENSICAEN, UNICAEN, GREYC, France) and Laurent Jacques (UCLouvain, Belgium).
%--
%main_09.tex
\item \hyperlink{main_09.1}{\large
``Low Rank and Group-Average Sparsity Driven Convex Optimization\vspace{1mm}\newline for Direct Exoplanets Imaging''.}\dotfill{\scriptsize (p.~\pageref*{pdf:main_09})}
\vspace{2mm} \newline
Beno\^{i}t Pairet (UCLouvain, Belgium),
    Laurent Jacques (UCLouvain, Belgium), Carlos A. Gomez Gonzalez
    (ULg, Belgium), Olivier Absil (ULg, Belgium).
%--
%main_10.tex
\item \hyperlink{main_10.1}{\large
``A fast algorithm for high-order sparse linear prediction''.}\dotfill{\scriptsize (p.~\pageref*{pdf:main_10})}
\vspace{2mm} \newline
Tobias Lindstr\o m Jensen (Aalborg
    Universitet, Denmark), Daniele Giacobello (DTS Inc., Calabasas,
    CA, USA), Toon van Waterschoot (KULeuven, Belgium), Mads Gr\ae
    sb\o ll Christensen (Aalborg Universitet, Denmark).
%--
%main_11.tex
\item \hyperlink{main_11.1}{\large
``Compressive Hyperspectral Imaging with Fourier Transform Interferometry''.}\dotfill{\scriptsize (p.~\pageref*{pdf:main_11})}
\vspace{2mm} \newline
A. Moshtaghpour (UCLouvain, Belgium),
    K. Degraux (UCLouvain, Belgium), V. Cambareri (UCLouvain,
    Belgium), A. Gonzalez (UCLouvain, Belgium), M. Roblin (Lambda-X,
    Belgium), L. Jacques (UCLouvain, Belgium),  and P. Antoine (Lambda-X, Belgium).
%--
%main_13.tex
\item \hyperlink{main_13.1}{\large
``Inferring Sparsity: Compressed Sensing using Generalized\vspace{1mm}\newline Restricted Boltzmann Machines''.}\dotfill{\scriptsize (p.~\pageref*{pdf:main_13})}
\vspace{2mm} \newline
Eric W. Tramel (Ecole Normale Sup\'erieure,
    PSL Research University, France),
                          Andre Manoel (Ecole Normale Sup\'erieure,
    PSL Research University, France),
                          Francesco Caltagirone (INRIA Paris), 
                          Marylou Gabri\'e (Ecole Normale Sup\'erieure,
    PSL Research University, France) and
                          Florent Krzakala (Ecole Normale Sup\'erieure,
    PSL Research University, France).
%--
%main_14.tex
\item \hyperlink{main_14.1}{\large
``Interpolation on manifolds using B\'ezier functions''.}\dotfill{\scriptsize (p.~\pageref*{pdf:main_14})}
\vspace{2mm} \newline
Pierre-Yves Gousenbourger (UCLouvain, Belgium), P.-A. Absil
(UCLouvain, Belgium), Benedikt Wirth (U. M\"unster, Germany) and Laurent Jacques (UCLouvain, Belgium).
%--
%main_15.tex
\item \hyperlink{main_15.1}{\large
``Reliable recovery of hierarchically sparse signals''.}\dotfill{\scriptsize (p.~\pageref*{pdf:main_15})}
\vspace{2mm} \newline
Ingo Roth (Freie Universit\"at Berlin, Germany), 
Martin Kliesch (Freie Universit\"at Berlin, Germany), 
Gerhard Wunder (Freie Universit\"at Berlin, Germany),
and Jens Eisert (Freie Universit\"at Berlin, Germany).
%--
%main_16.tex
\item \hyperlink{main_16.1}{\large
``Minimizing Isotropic Total Variation without Subiterations''.}\dotfill{\scriptsize (p.~\pageref*{pdf:main_16})}
\vspace{2mm} \newline
Ulugbek~S.~Kamilov (MERL, USA).
%--
%main_17.tex
\item \hyperlink{main_17.1}{\large
``Learning MMSE Optimal Thresholds for FISTA''.}\dotfill{\scriptsize (p.~\pageref*{pdf:main_17})}
\vspace{2mm} \newline
Ulugbek~S.~Kamilov (MERL, USA) and Hassan Mansour (MERL, USA).
%--
%main_18.tex
\item \hyperlink{main_18.1}{\large
``The best of both worlds: synthesis-based acceleration\vspace{1mm}\newline for physics-driven cosparse regularization''.}\dotfill{\scriptsize (p.~\pageref*{pdf:main_18})}
\vspace{2mm} \newline
Sr\dj{}an Kiti\'c (Technicolor R\&D, France),
    Nancy Bertin (CNRS - UMR 6074, France) and R\'emi Gribonval (Inria, France).
%--
%main_19.tex
\item \hyperlink{main_19.1}{\large
``A Student-t based sparsity enforcing  hierarchical prior for linear inverse problems\vspace{1mm}\newline and its efficient Bayesian computation for 2D and 3D Computed Tomography''.}\dotfill{\scriptsize (p.~\pageref*{pdf:main_19})}
\vspace{2mm} \newline
Ali Mohammad-Djafari
    (CentraleSup\'elec-U. Paris Saclay, France), Li Wang
    (CentraleSup\'elec-U. Paris Saclay, France), Nicolas Gac
    (CentraleSup\'elec-U. Paris Saclay, France) and Folkert Bleichrodt
    (CWI, The Netherlands). 
%--
%main_20.tex
\item \hyperlink{main_20.1}{\large
``Simultaneous reconstruction and separation in a spectral CT framework''.}\dotfill{\scriptsize (p.~\pageref*{pdf:main_20})}
\vspace{2mm} \newline
S. Tairi (CPPM, France), S. Anthoine (Aix
    Marseille Universit\'e, France), C. Morel (CPPM, France) and
    Y. Boursier (CPPM, France).
%--
%main_21.tex
\item \hyperlink{main_21.1}{\large
``Debiasing incorporated into reconstruction of low-rank modelled dynamic MRI data''.}\dotfill{\scriptsize (p.~\pageref*{pdf:main_21})}
\vspace{2mm} \newline
Marie Da\v{n}kov\'a (Brno University of Technology
    \& Masaryk University, Czech Republic) and Pavel Rajmic (Brno University of Technology, Czech Republic).
%--
%main_22.tex
\item \hyperlink{main_22.1}{\large
``Sparse MRI with a Markov Random Field Prior for the Subband Coefficients''.}\dotfill{\scriptsize (p.~\pageref*{pdf:main_22})}
\vspace{2mm} \newline
Marko Pani\'{c} (University of Novi Sad, Serbia),
Dejan Vukobratovic (University of Novi Sad, Serbia),
Vladimir Crnojevi\'{c} (University of Novi Sad, Serbia) and Aleksandra
Pi\v{z}urica (Ghent University, Belgium).
%--
%main_23.tex
\item \hyperlink{main_23.1}{\large
``Active GAP screening for the LASSO''.}\dotfill{\scriptsize (p.~\pageref*{pdf:main_23})}
\vspace{2mm} \newline
A. Bonnefoy (Aix Marseille Universit\'e, France) and S. Anthoine (Aix Marseille Universit\'e, France).
%--
%main_24.tex
\item \hyperlink{main_24.1}{\large
``Paint Loss Detection in Old Paintings by Sparse Representation
Classification''.}\dotfill{\scriptsize (p.~\pageref*{pdf:main_24})} 
\vspace{2mm} \newline
Shaoguang Huang (Ghent University, Belgium),
Wenzhi Liao (Ghent University, Belgium), Hongyan Zhang (Wuhan
University, China) and Aleksandra Pi\v{z}urica (Ghent University, Belgium).
%%%%%%%%%%%%%%%%%%%%%%%%%
\end{itemize}

\newpage
\setcounter{page}{1}
\renewcommand{\thepage}{\arabic{page}}

\label{pdf:main_02}
\includepdf[offset=.65cm -2mm,pages=-,link=true,linkname=main_02,pagecommand={}]{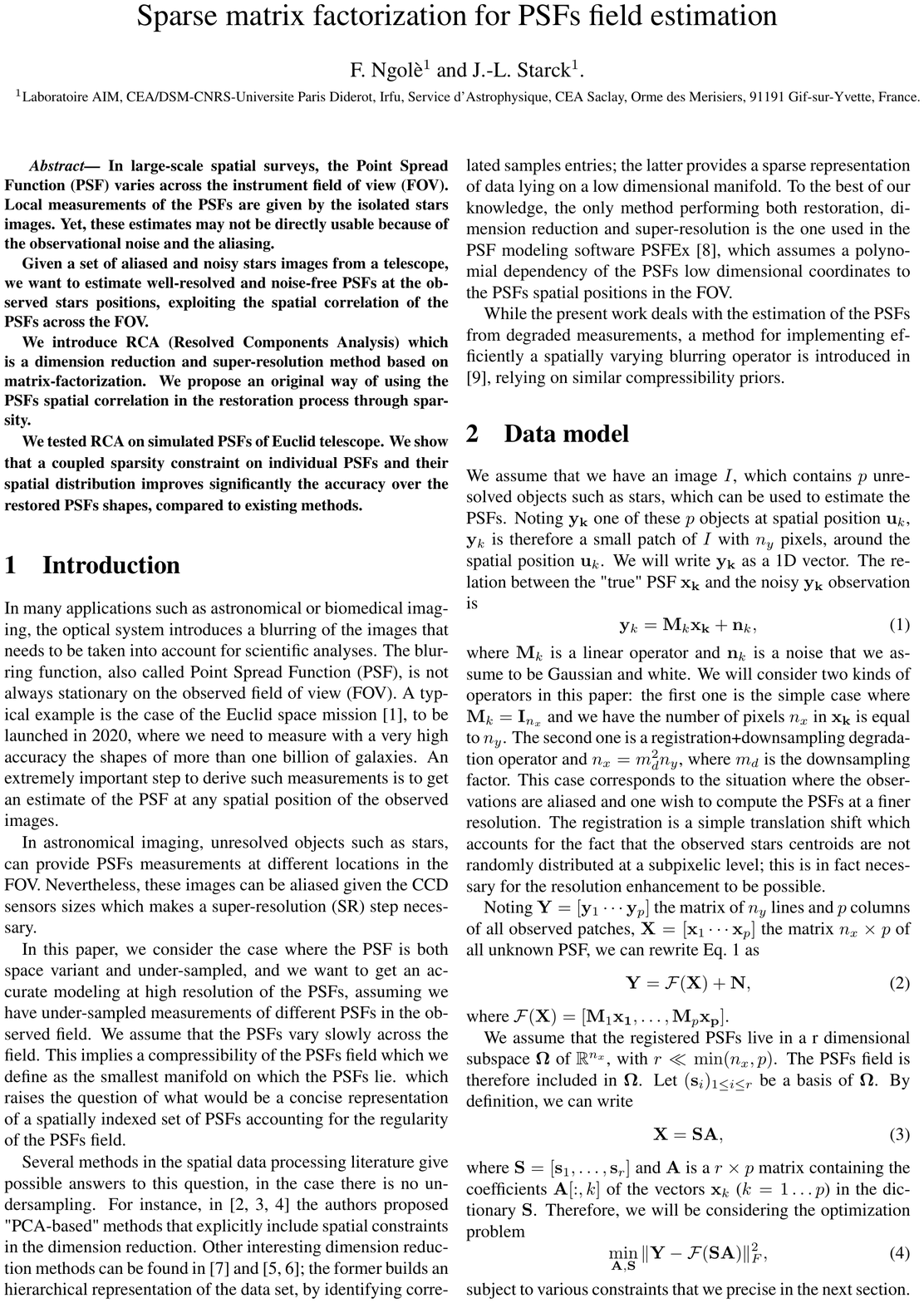}

\label{pdf:main_03}
\includepdf[offset=1cm -2mm,pages=-,link=true,linkname=main_03,pagecommand={}]{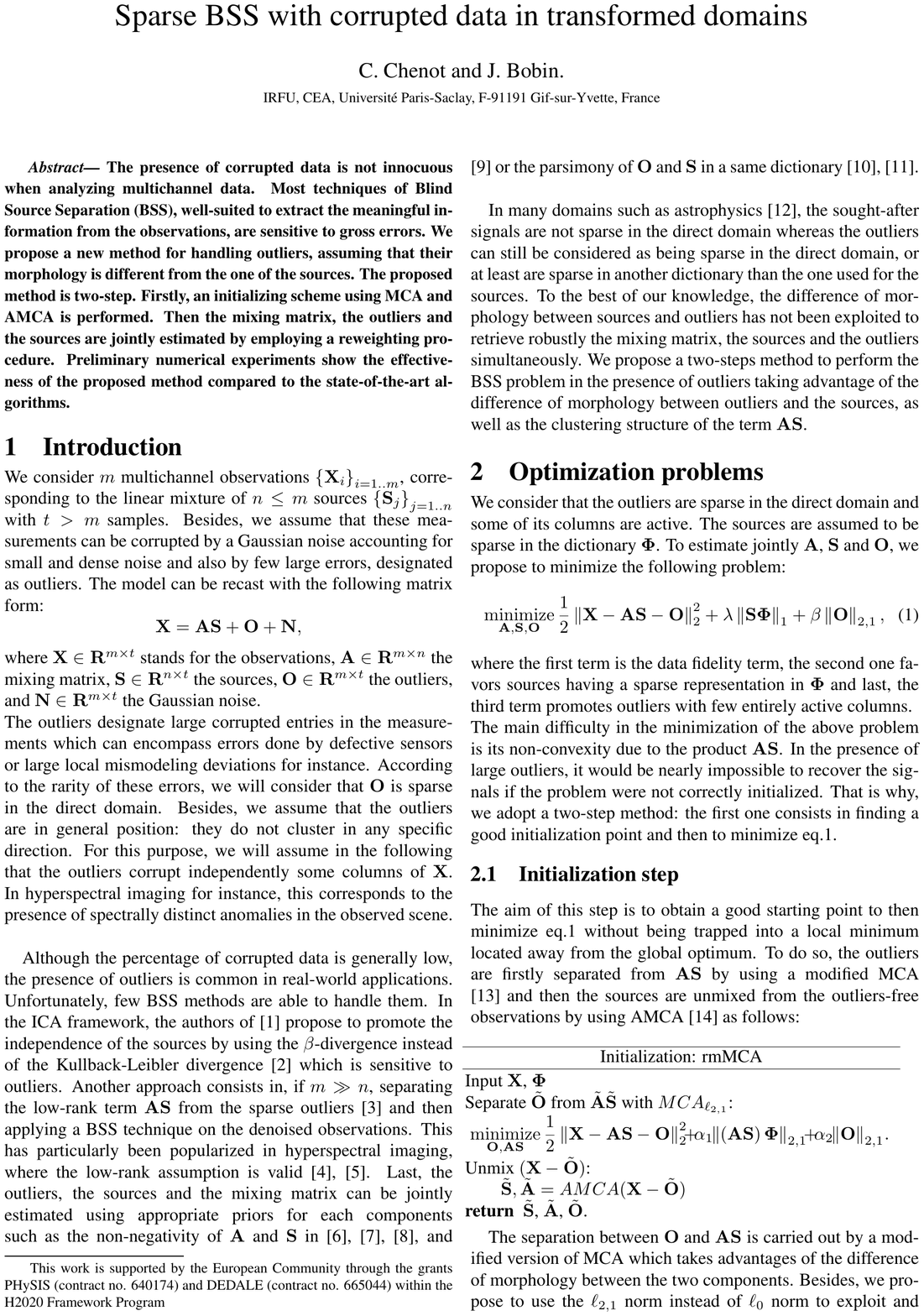}

\label{pdf:main_04}
\includepdf[offset=.65cm -2mm,pages=-,link=true,linkname=main_04,pagecommand={}]{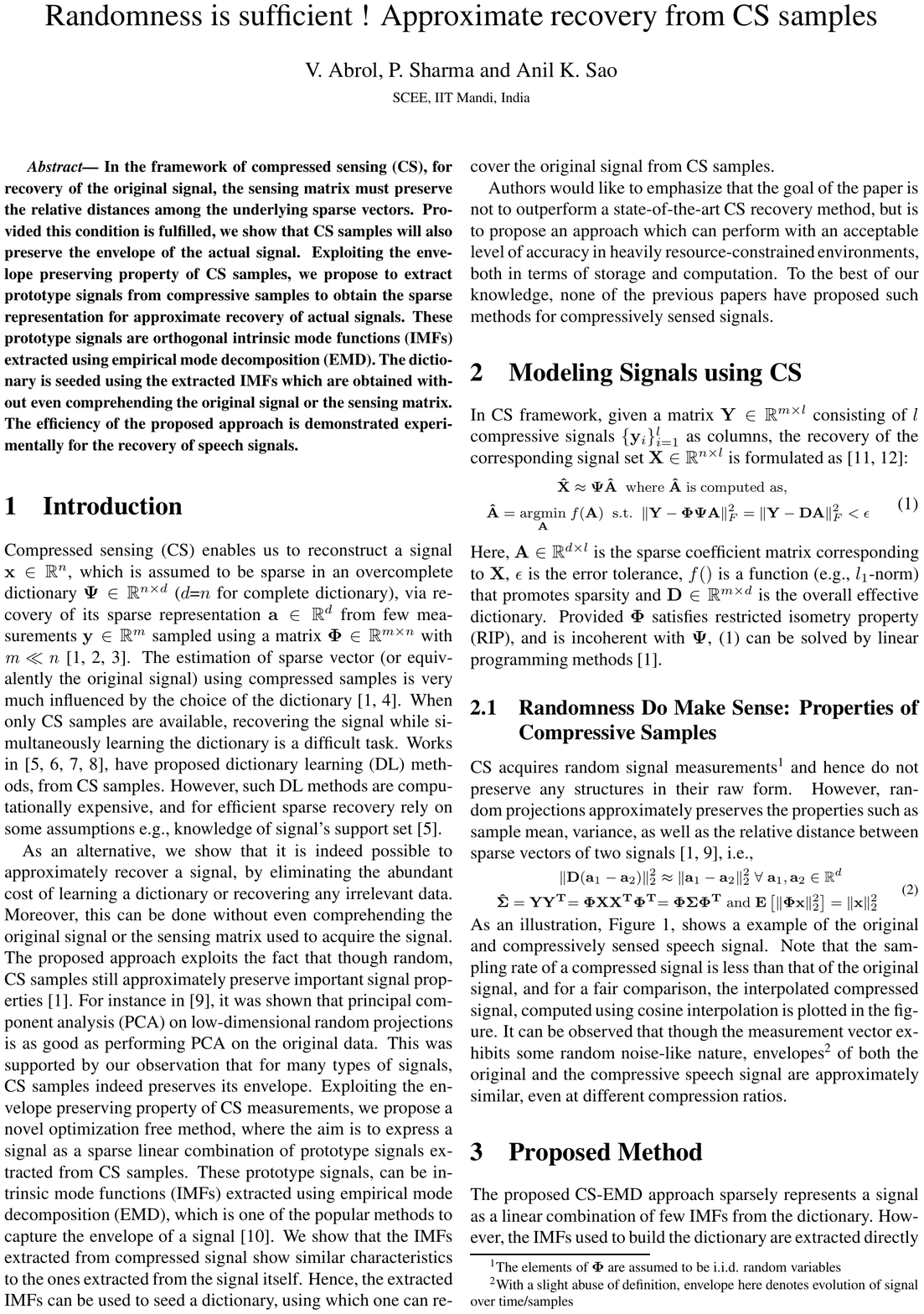}

\label{pdf:main_05}
\includepdf[offset=.65cm -2mm,pages=-,link=true,linkname=main_05,pagecommand={}]{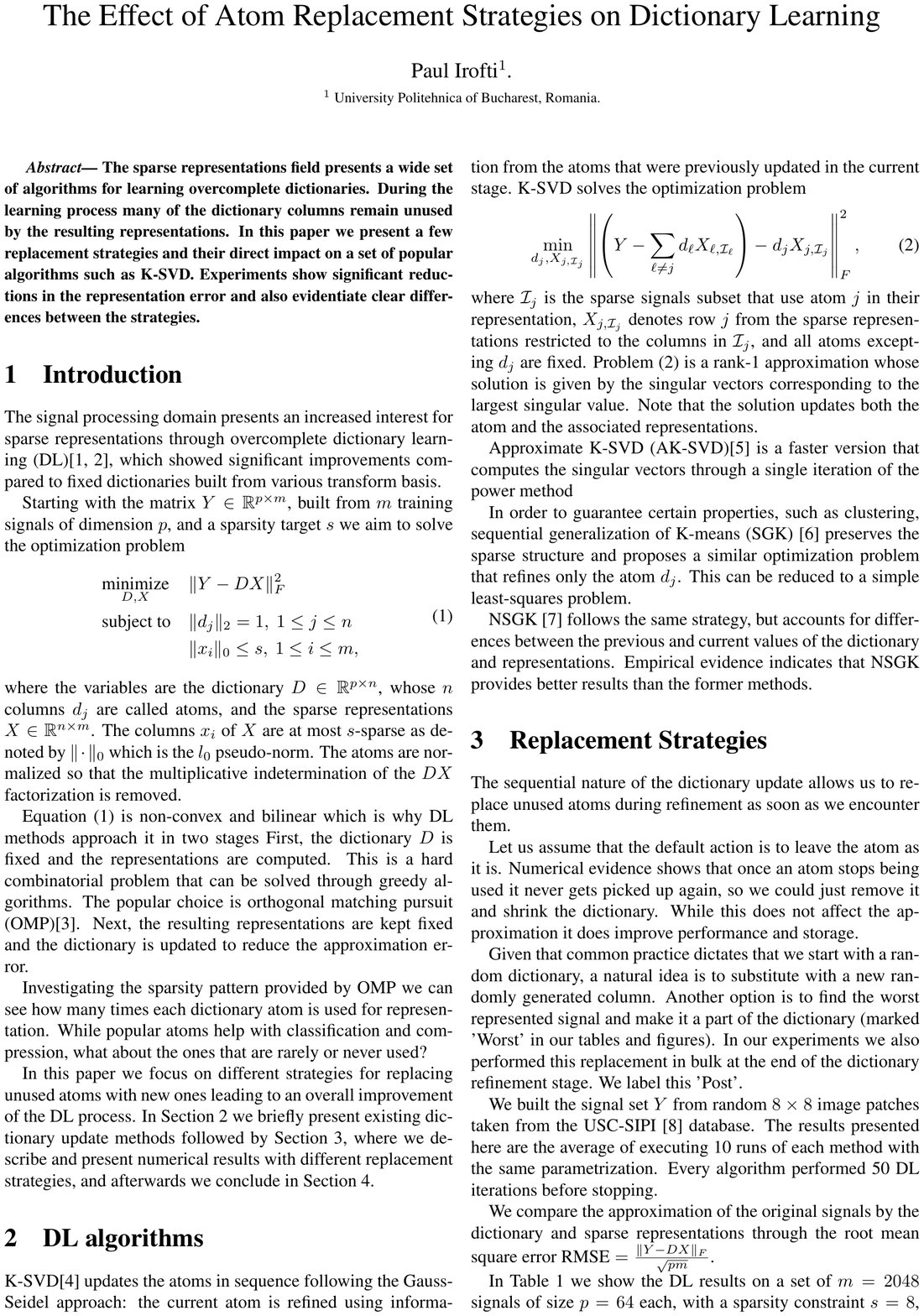}

\label{pdf:main_06}
\includepdf[offset=.65cm -2mm,pages=-,link=true,linkname=main_06,pagecommand={}]{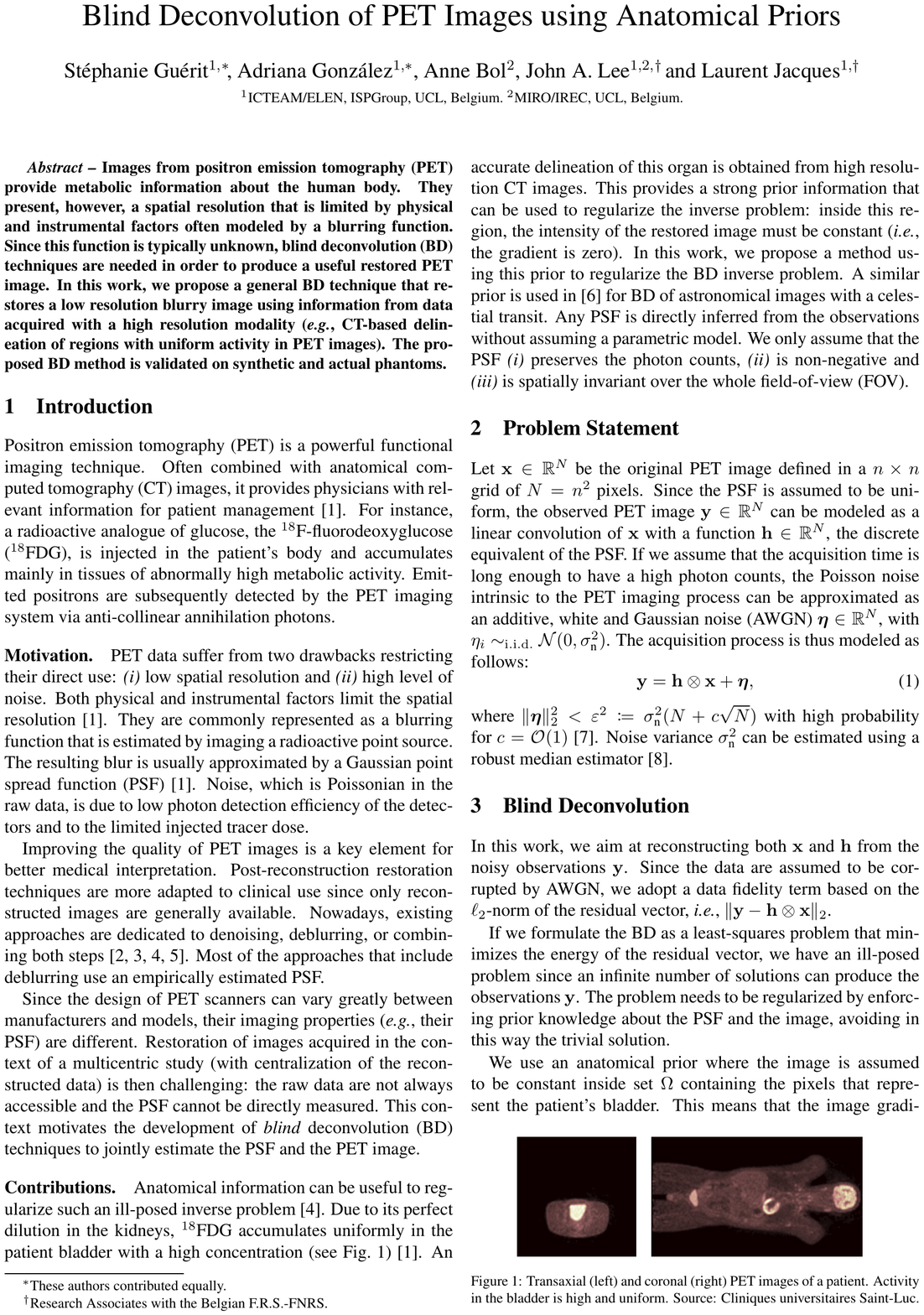}

\label{pdf:main_07}
\includepdf[offset=.65cm -2mm,pages=-,link=true,linkname=main_07,pagecommand={}]{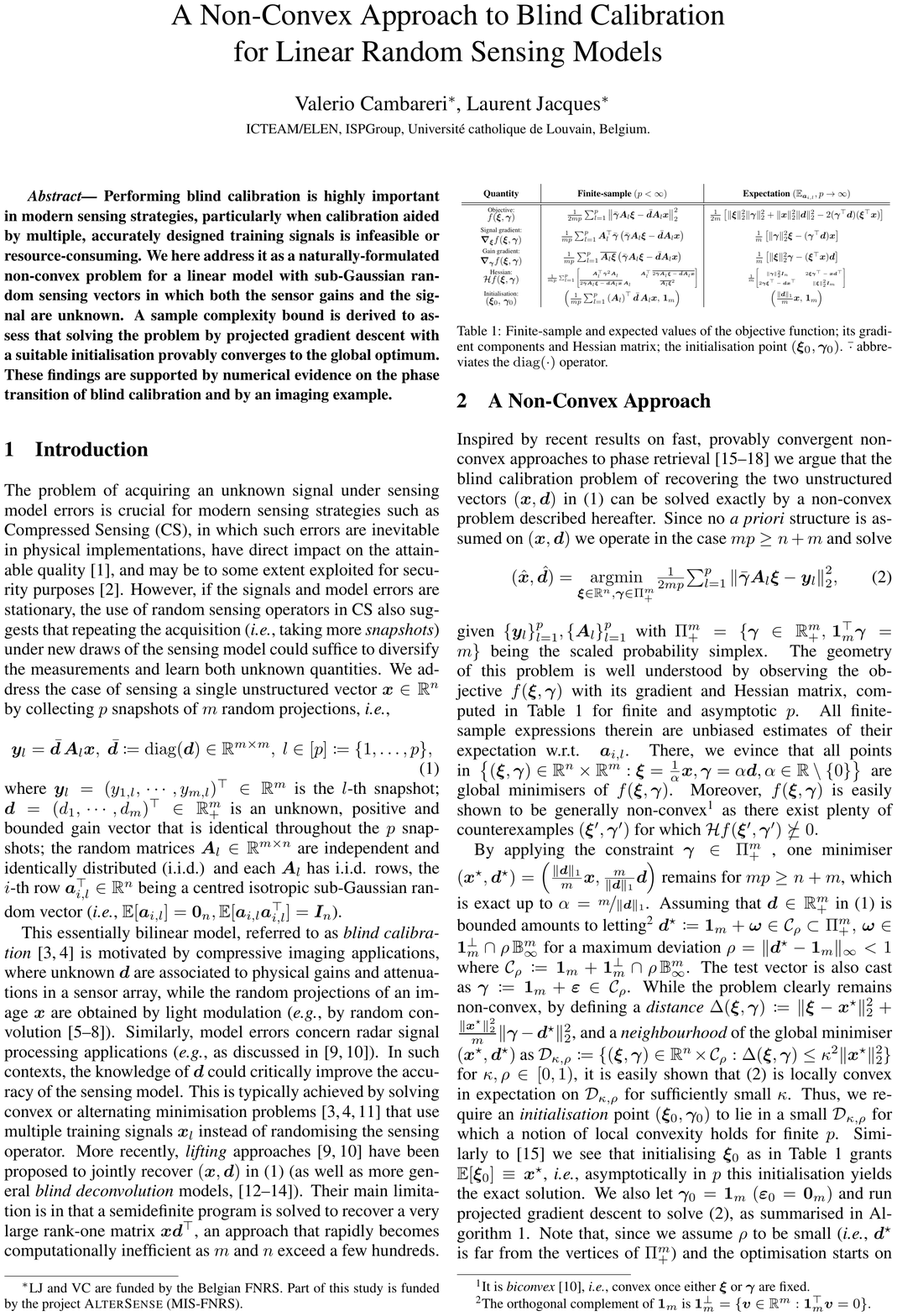}

\label{pdf:main_08}
\includepdf[offset=.65cm -2mm,pages=-,link=true,linkname=main_08,pagecommand={}]{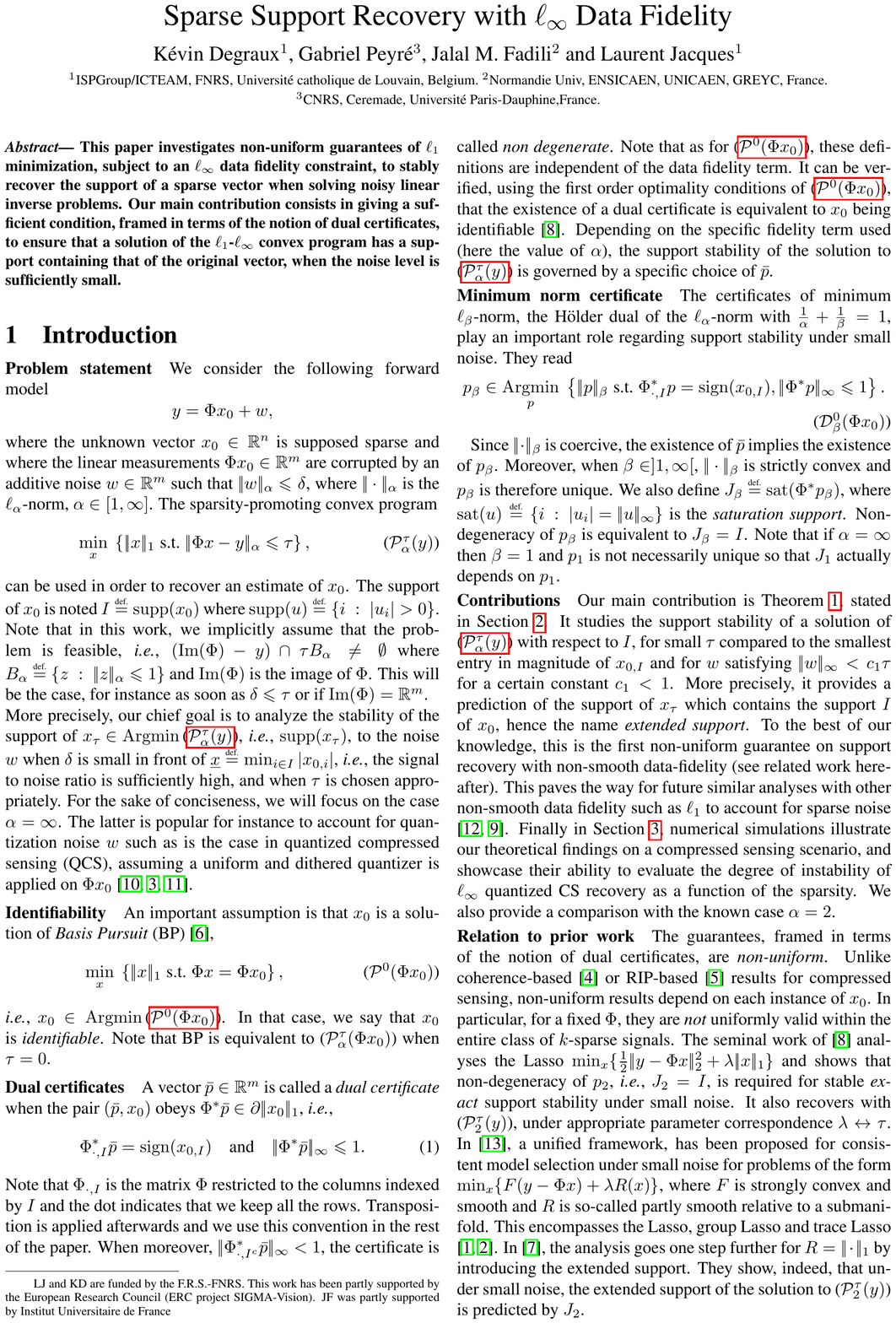}

\label{pdf:main_09}
\includepdf[offset=.65cm -2mm,pages=-,link=true,linkname=main_09,pagecommand={}]{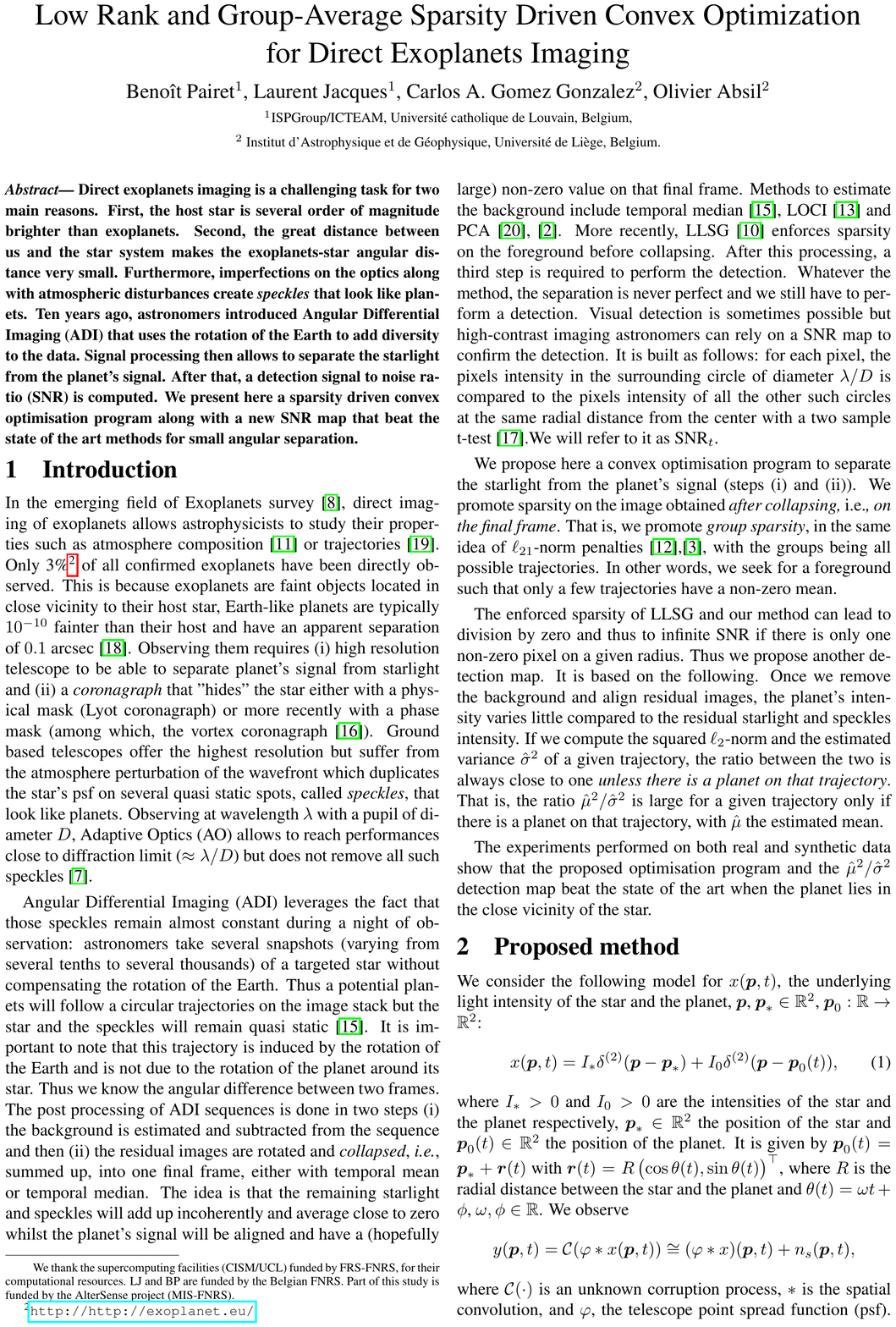}

\label{pdf:main_10}
\includepdf[offset=.65cm -2mm,pages=-,link=true,linkname=main_10,pagecommand={}]{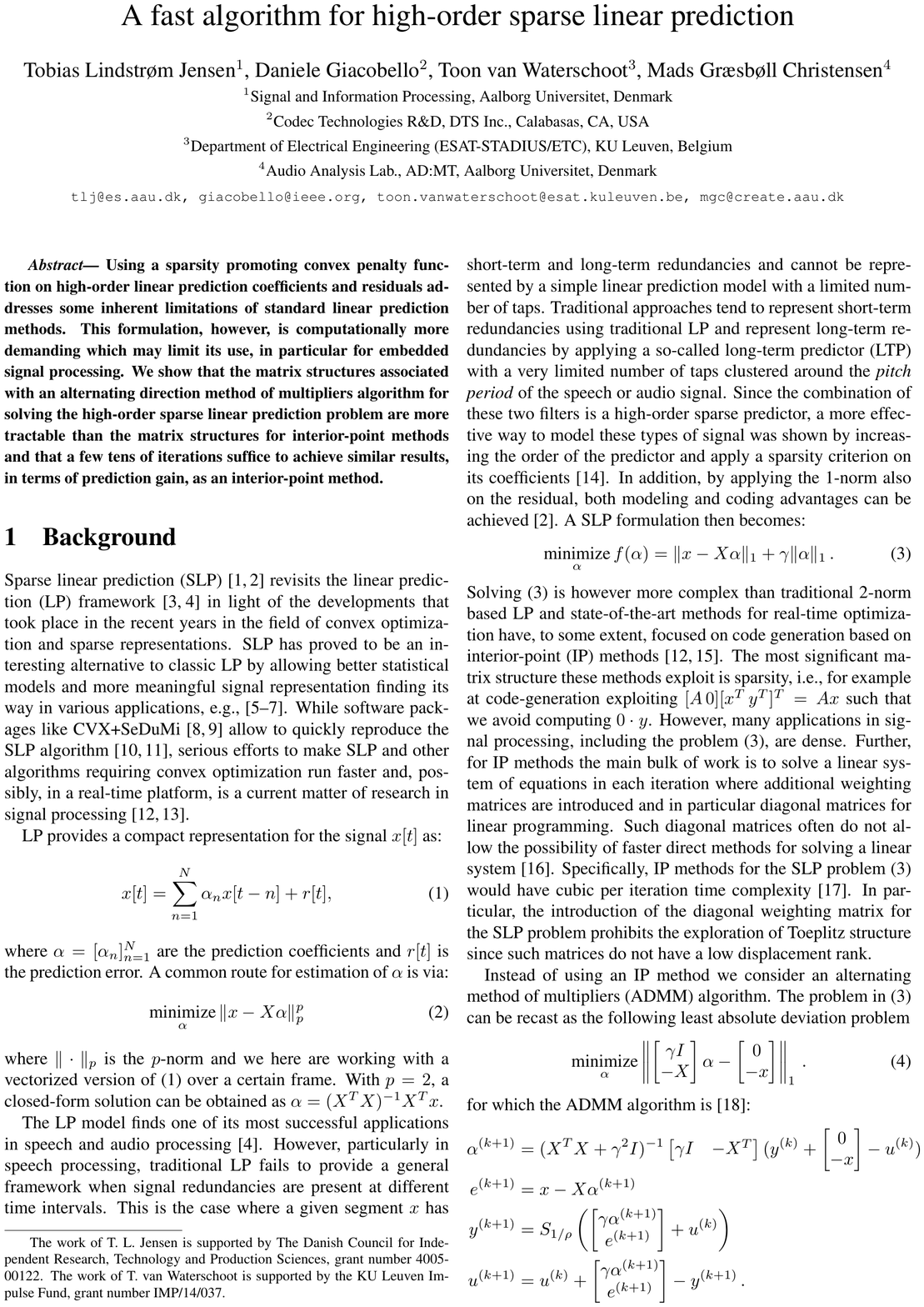}

\label{pdf:main_11}
\includepdf[offset=.65cm -2mm,pages=-,link=true,linkname=main_11,pagecommand={}]{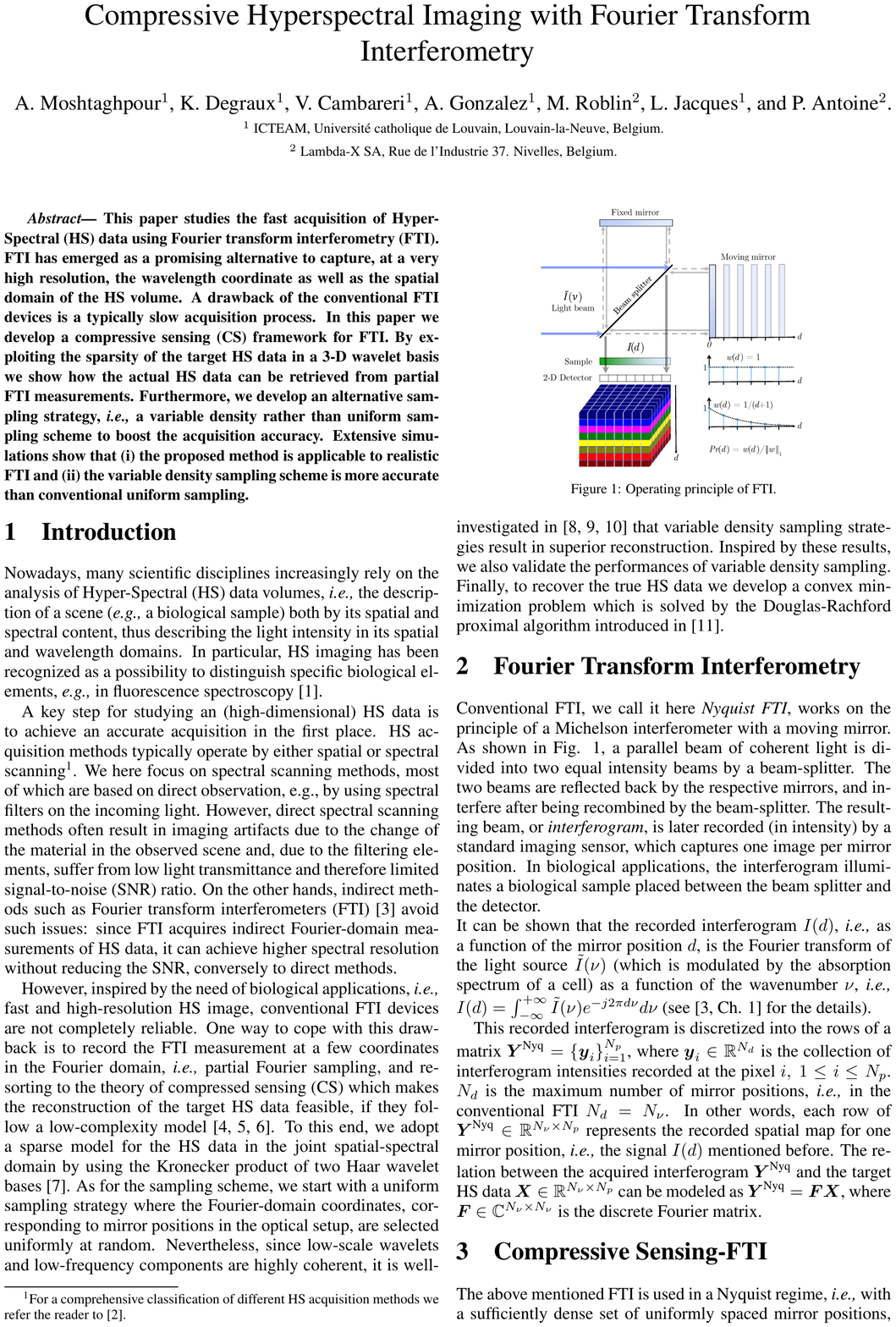}

\label{pdf:main_13}
\includepdf[offset=.65cm -2mm,pages=-,link=true,linkname=main_13,pagecommand={}]{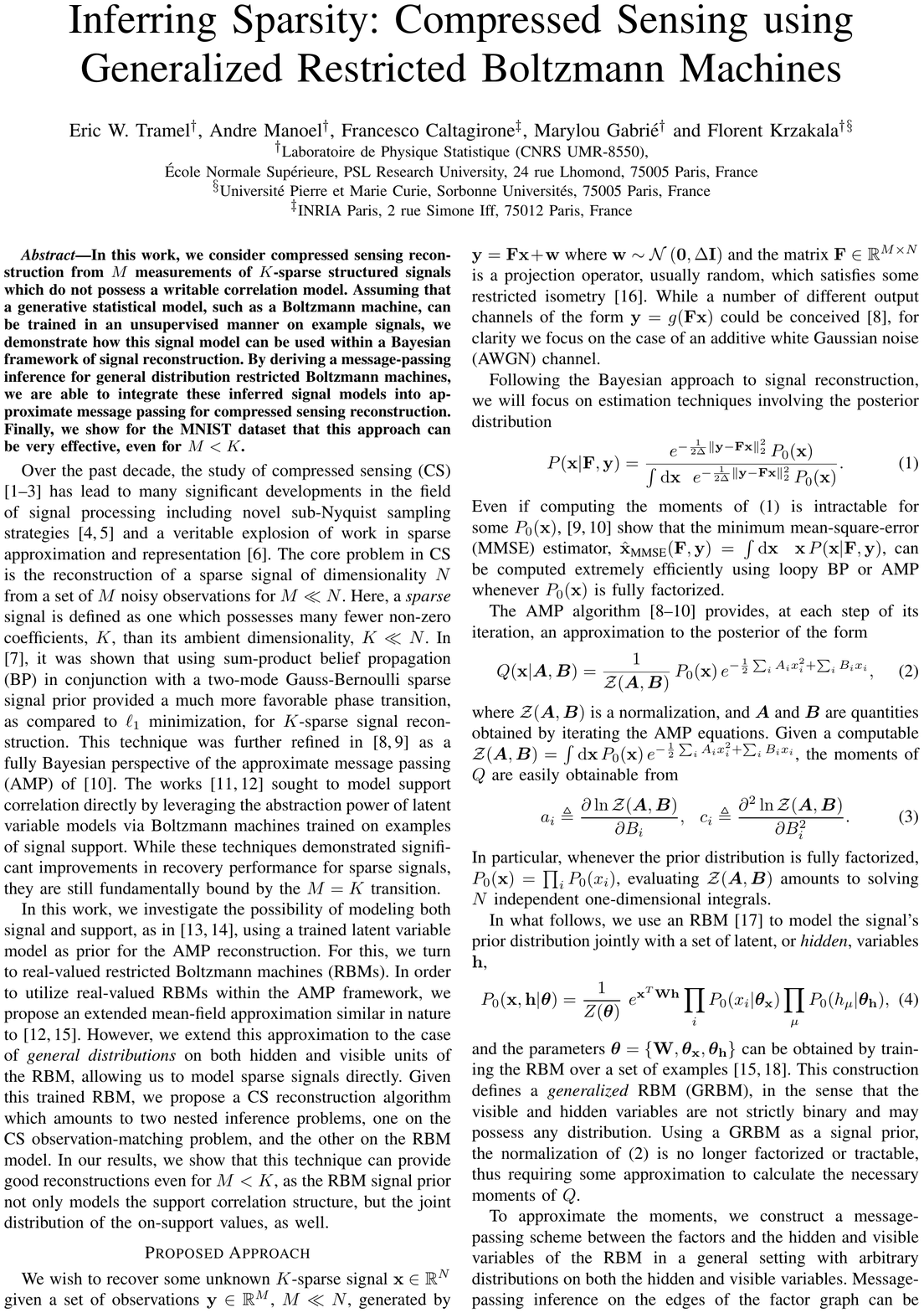}

\label{pdf:main_14}
\includepdf[offset=.65cm -2mm,pages=-,link=true,linkname=main_14,pagecommand={}]{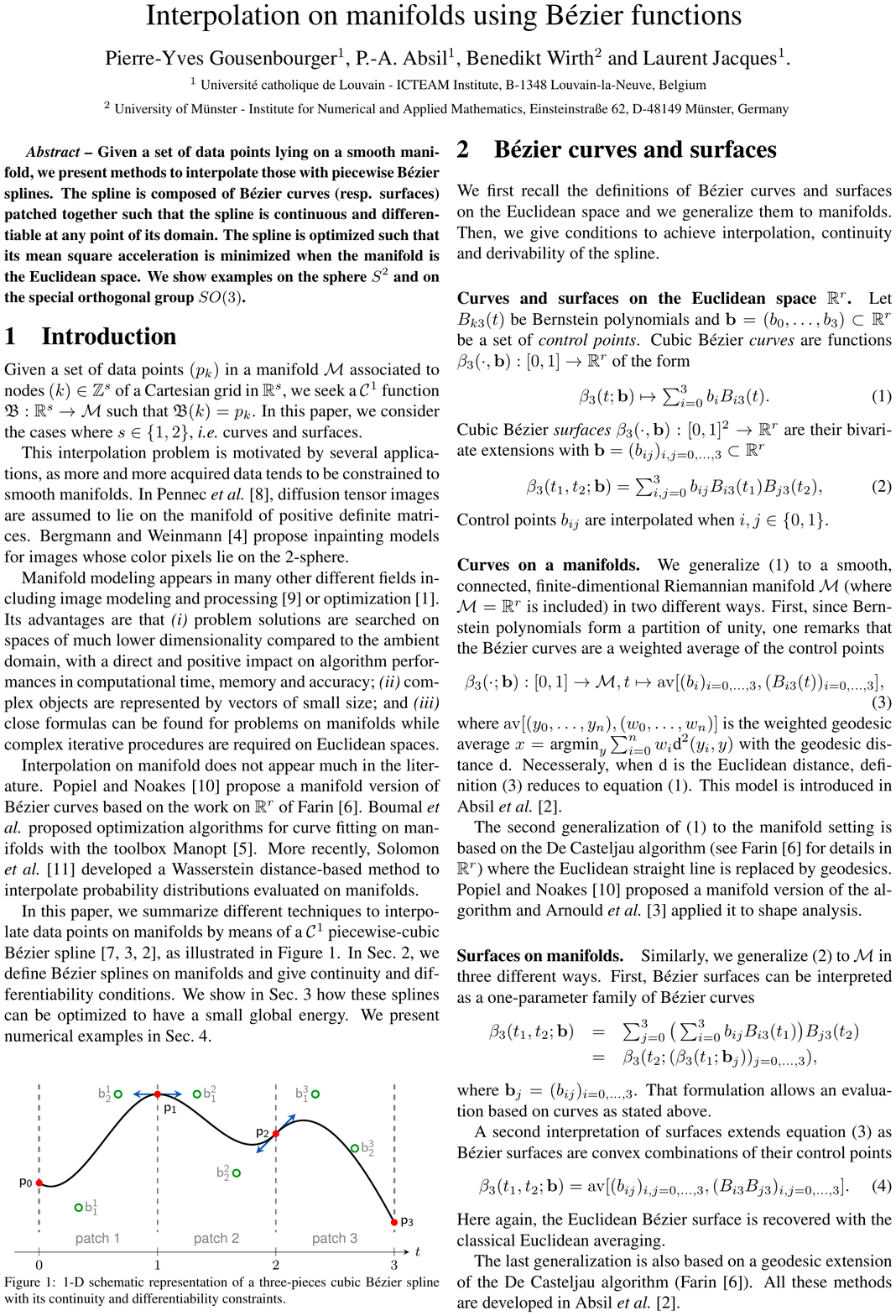}

\label{pdf:main_15}
\includepdf[offset=.65cm -2mm,pages=-,link=true,linkname=main_15,pagecommand={}]{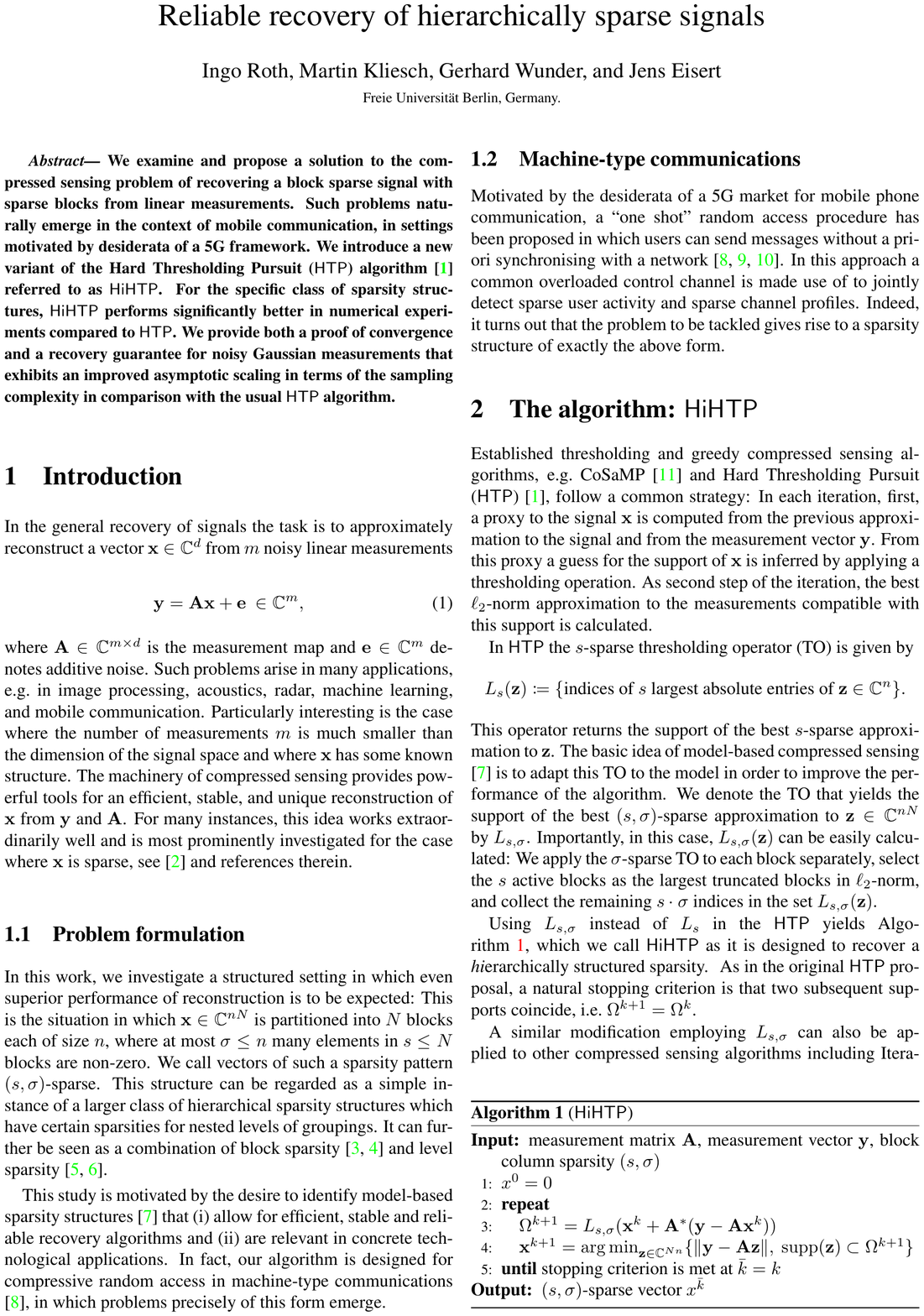}

\label{pdf:main_16}
\includepdf[offset=.65cm -2mm,pages=-,link=true,linkname=main_16,pagecommand={}]{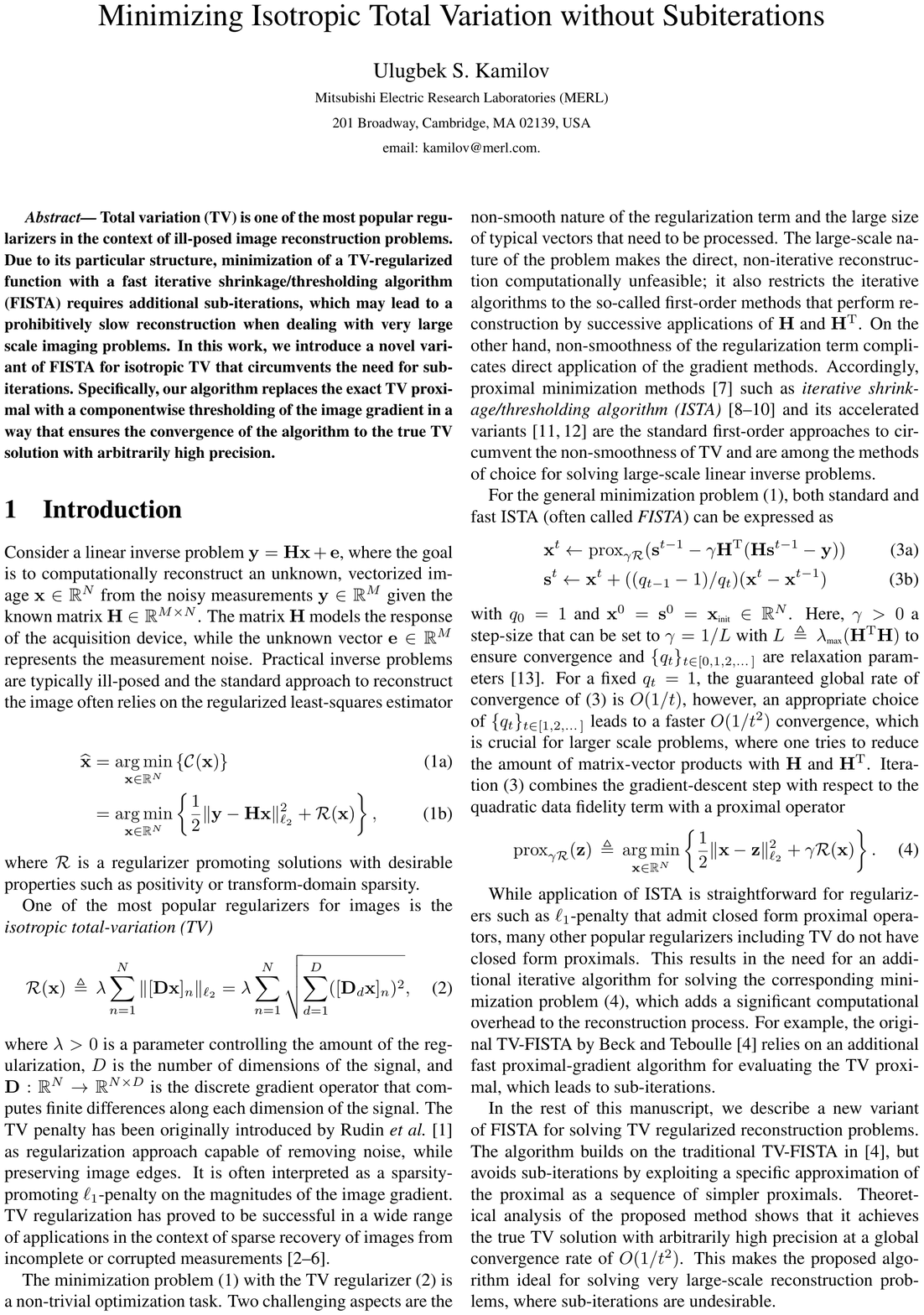}

\label{pdf:main_17}
\includepdf[offset=.65cm -2mm,pages=-,link=true,linkname=main_17,pagecommand={}]{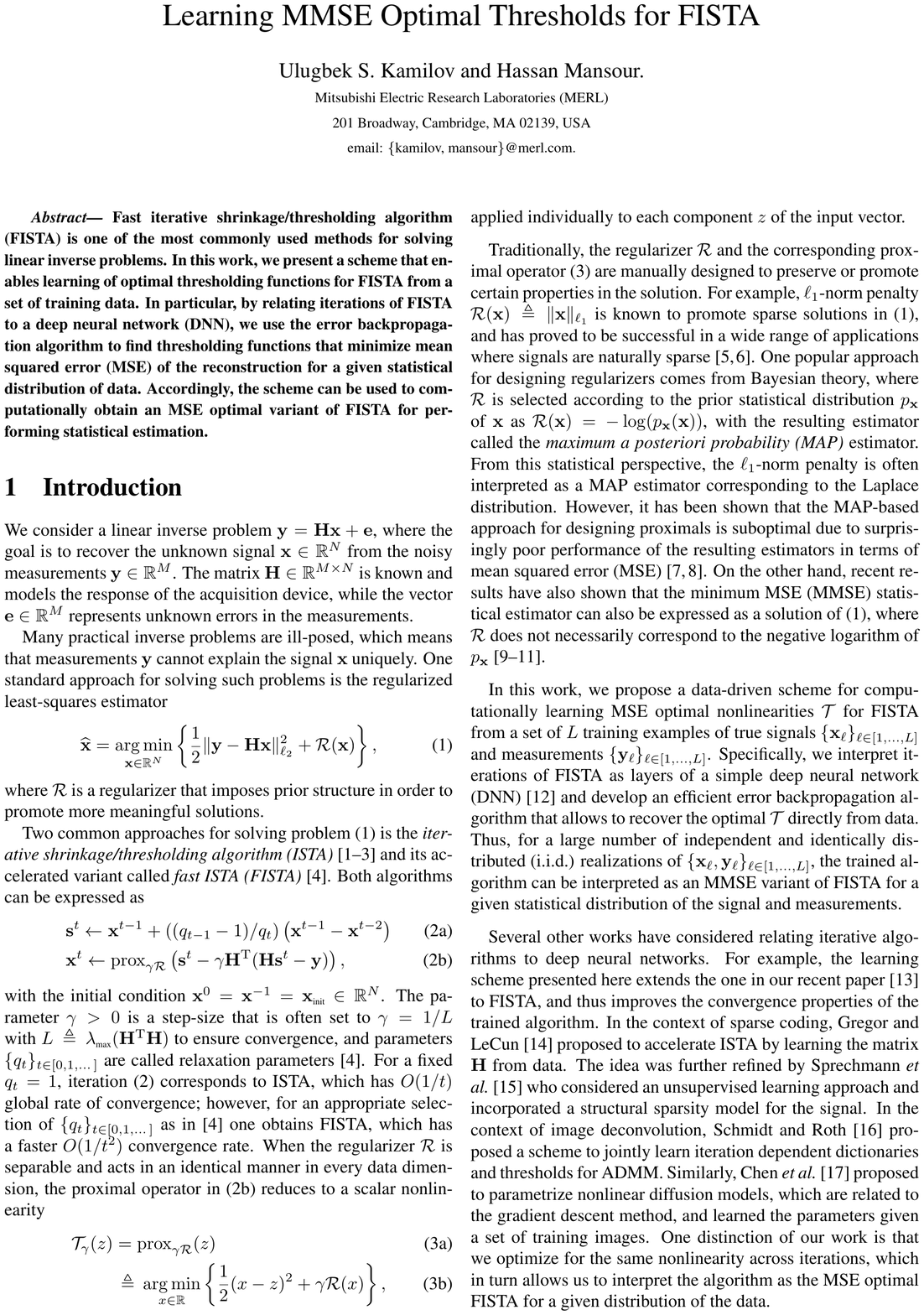}

\label{pdf:main_18}
\includepdf[offset=.65cm -2mm,pages=-,link=true,linkname=main_18,pagecommand={}]{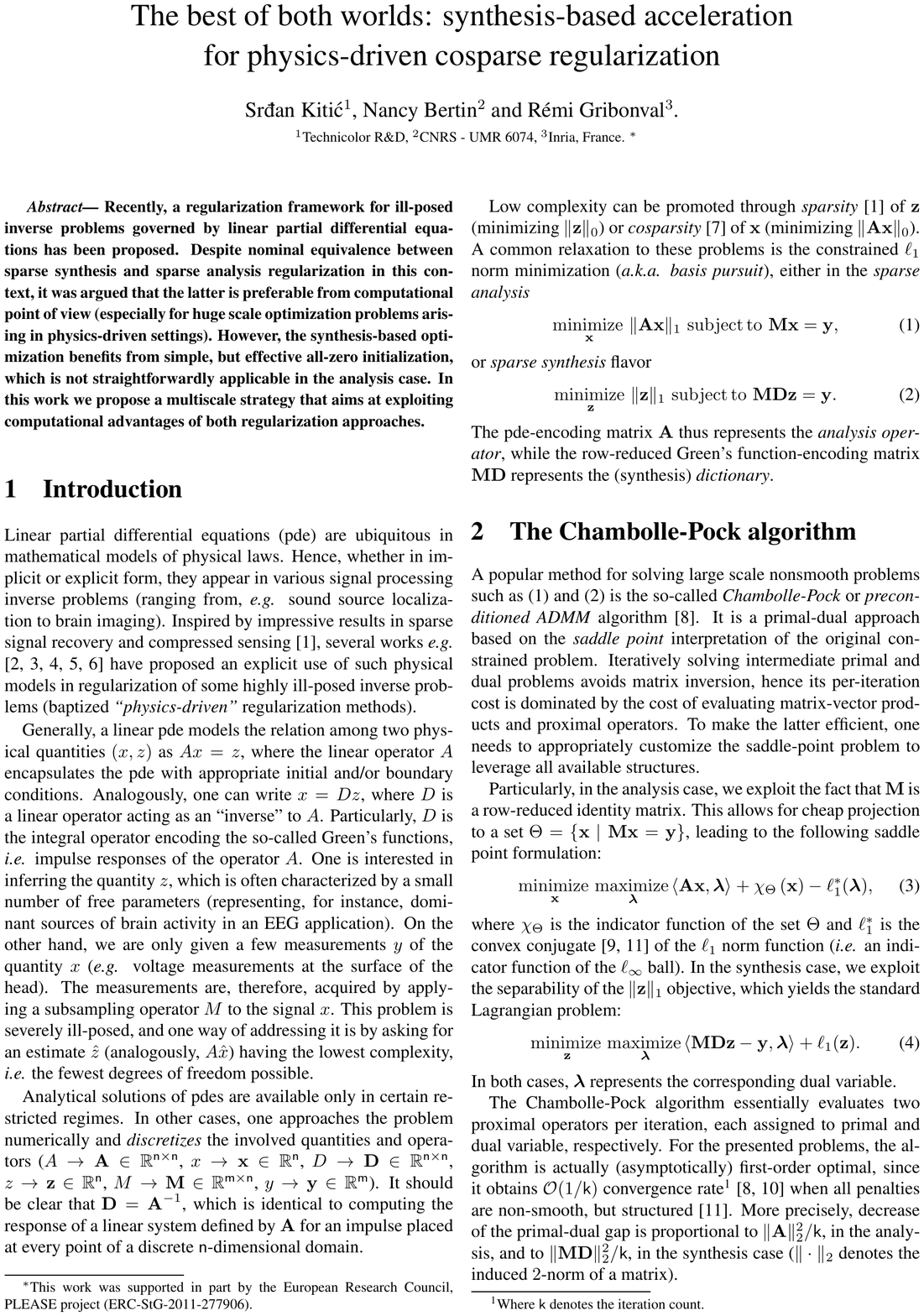}

\label{pdf:main_19}
\includepdf[offset=.65cm -2mm,pages=-,link=true,linkname=main_19,pagecommand={}]{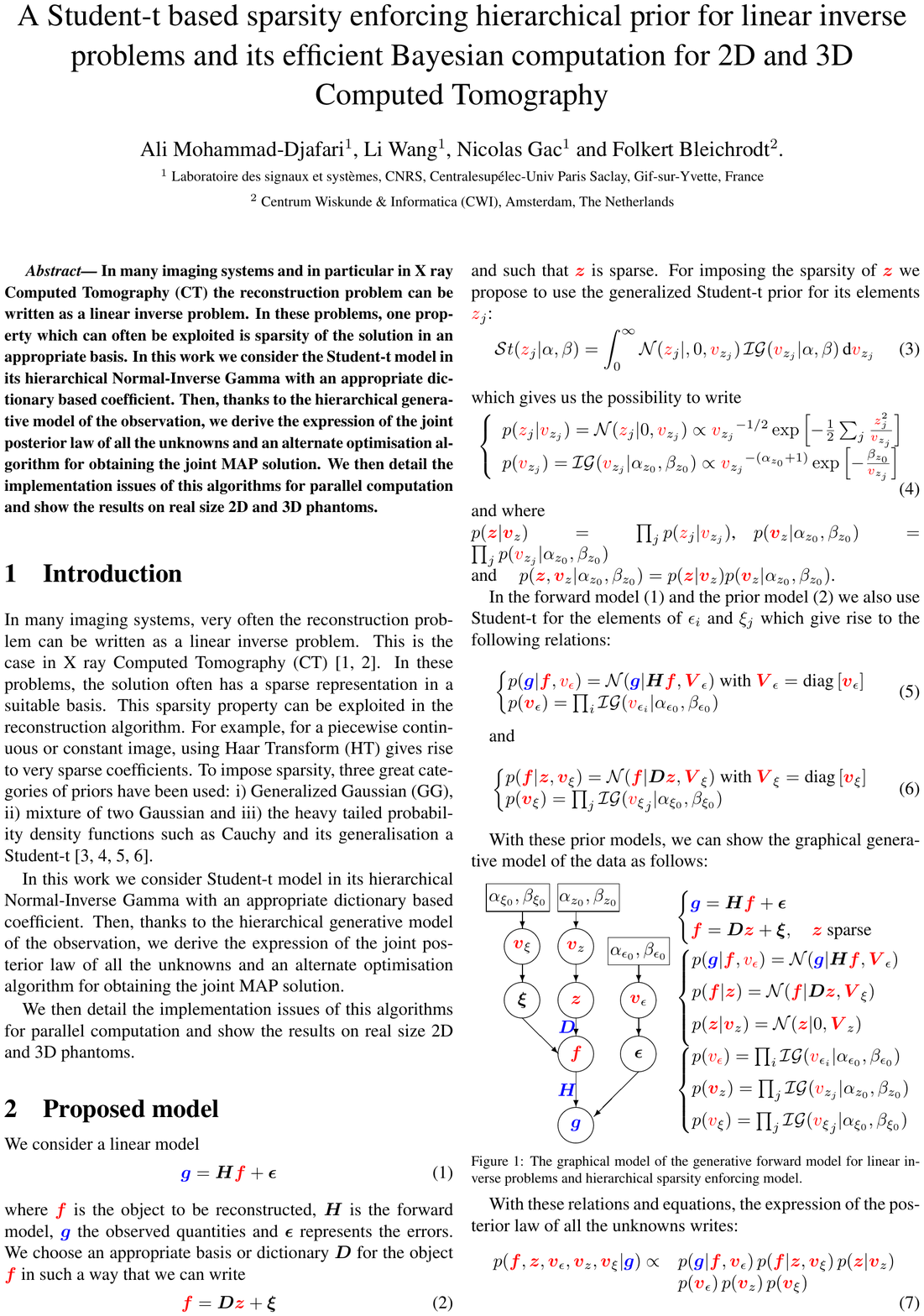}

\label{pdf:main_20}
\includepdf[offset=.65cm -2mm,pages=-,link=true,linkname=main_20,pagecommand={}]{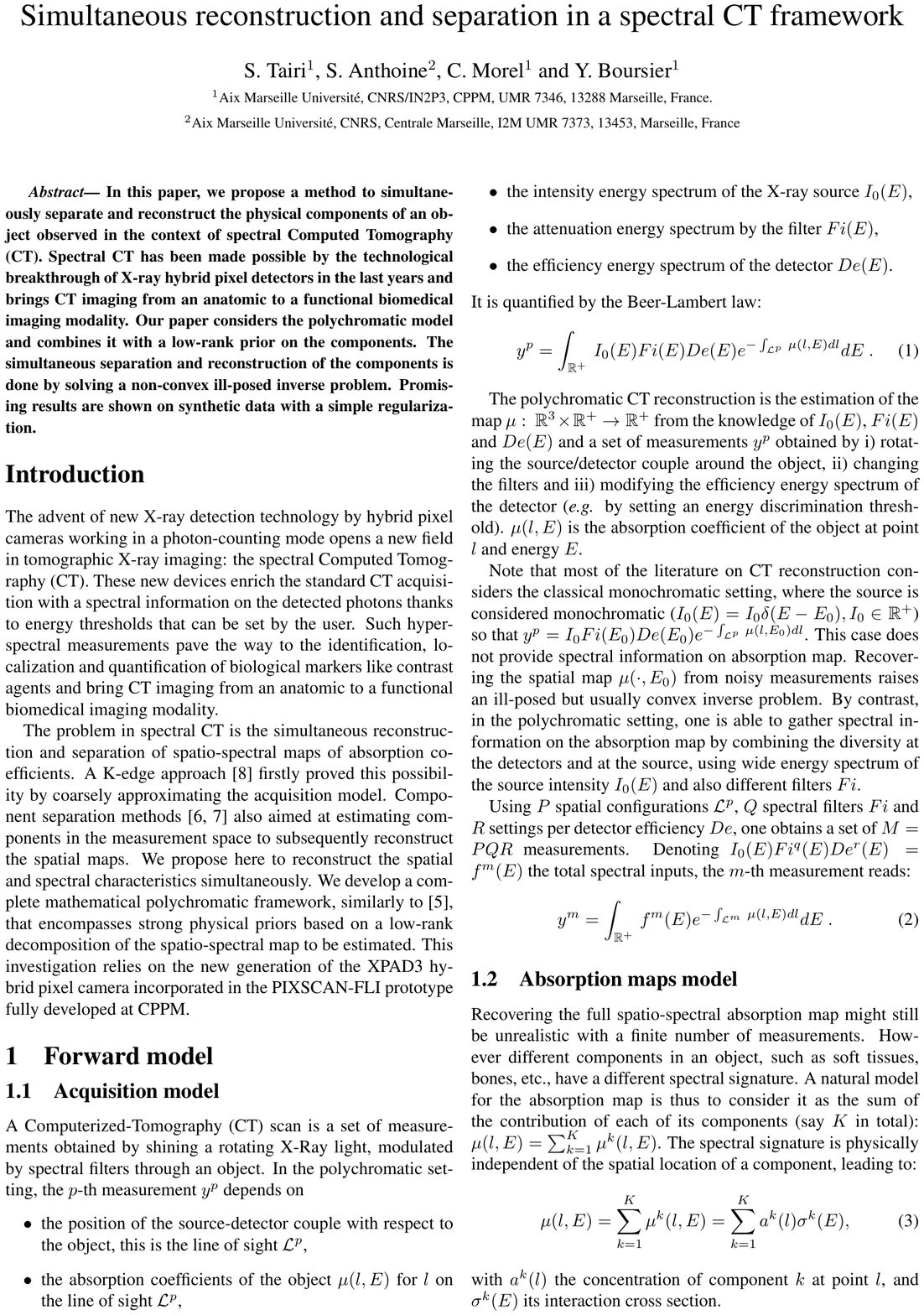}

\label{pdf:main_21}
\includepdf[offset=.65cm -2mm,pages=-,link=true,linkname=main_21,pagecommand={}]{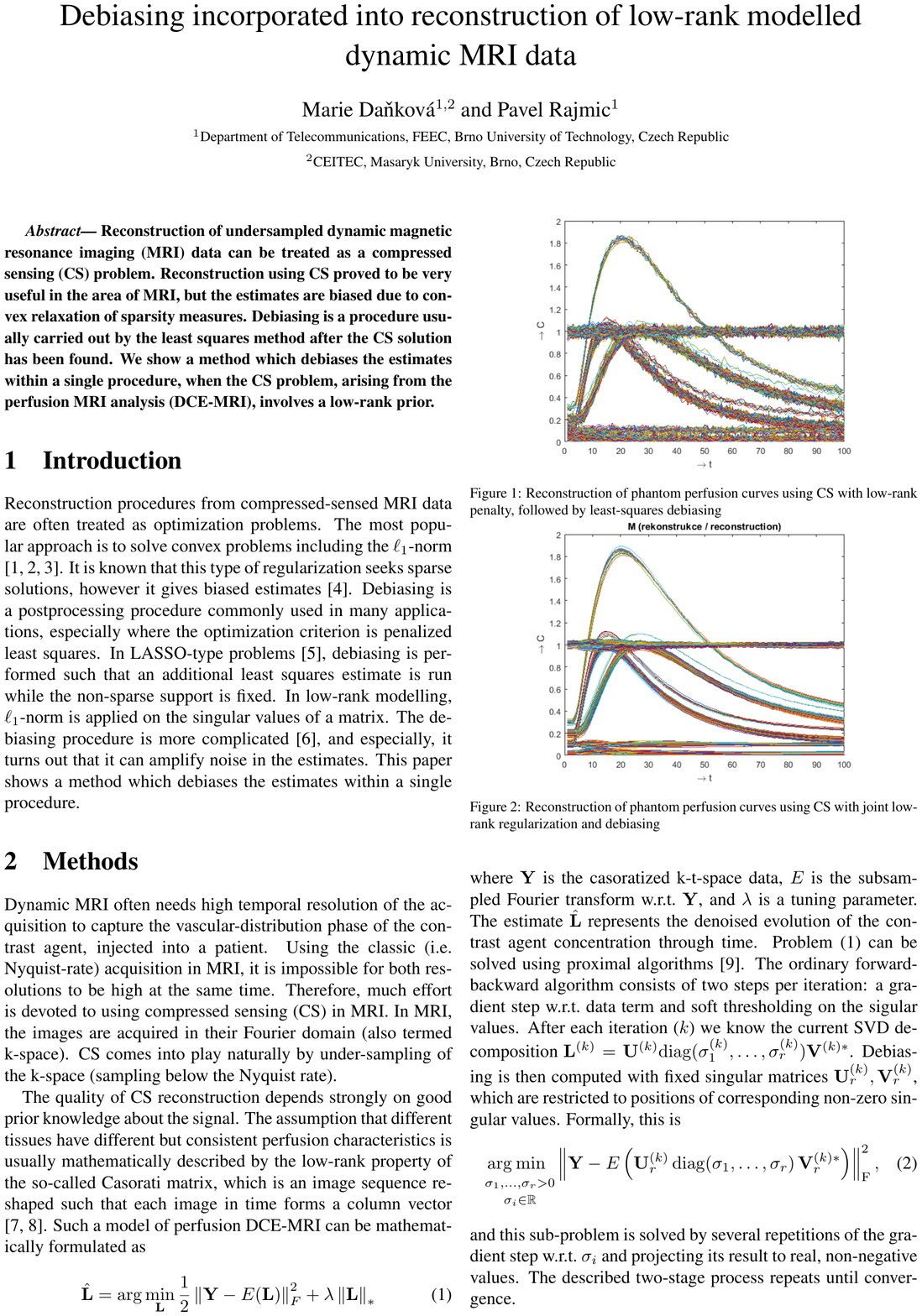}

\label{pdf:main_22}
\includepdf[offset=.65cm -2mm,pages=-,link=true,linkname=main_22,pagecommand={}]{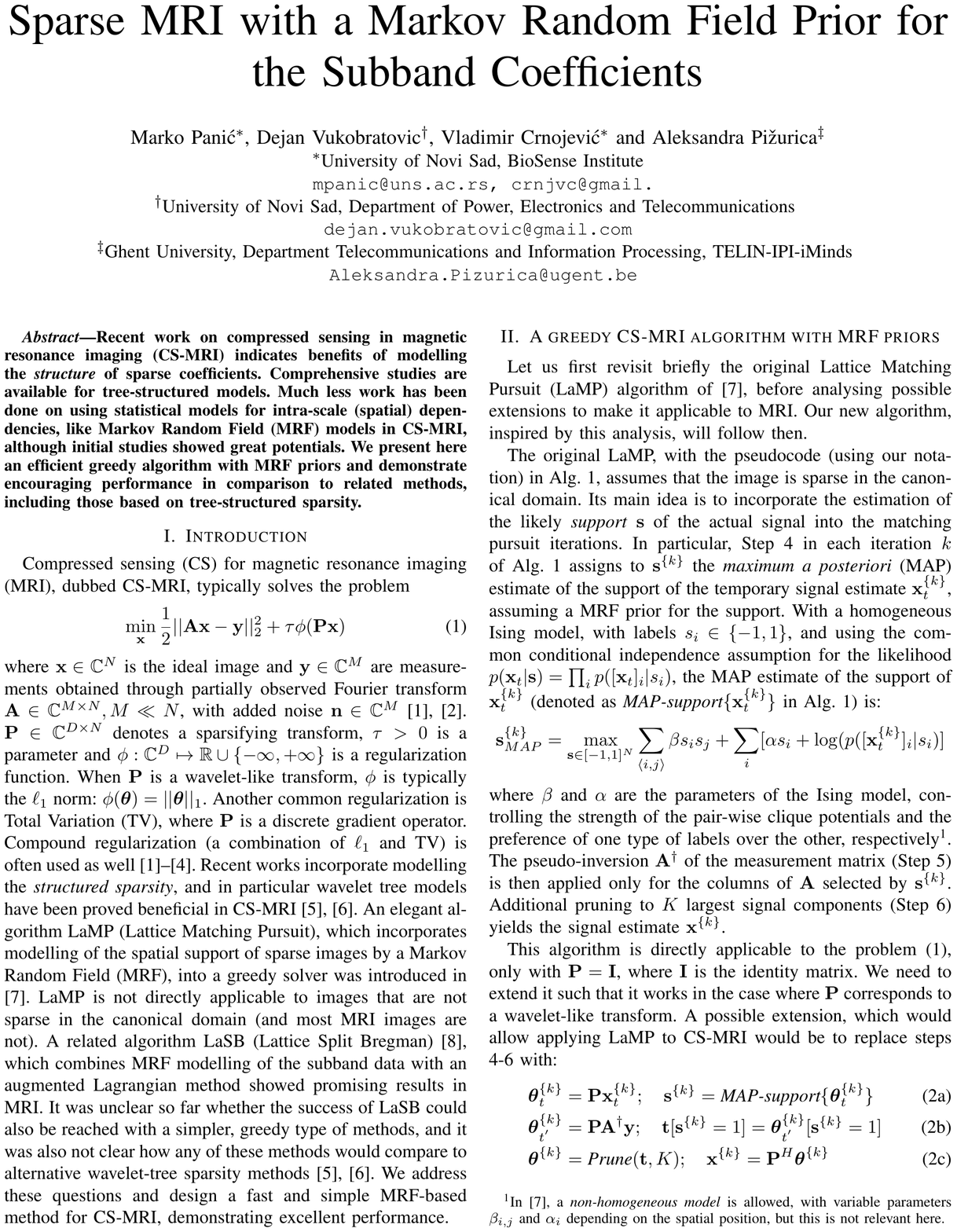}

\label{pdf:main_23}
\includepdf[offset=.65cm -2mm,pages=-,link=true,linkname=main_23,pagecommand={}]{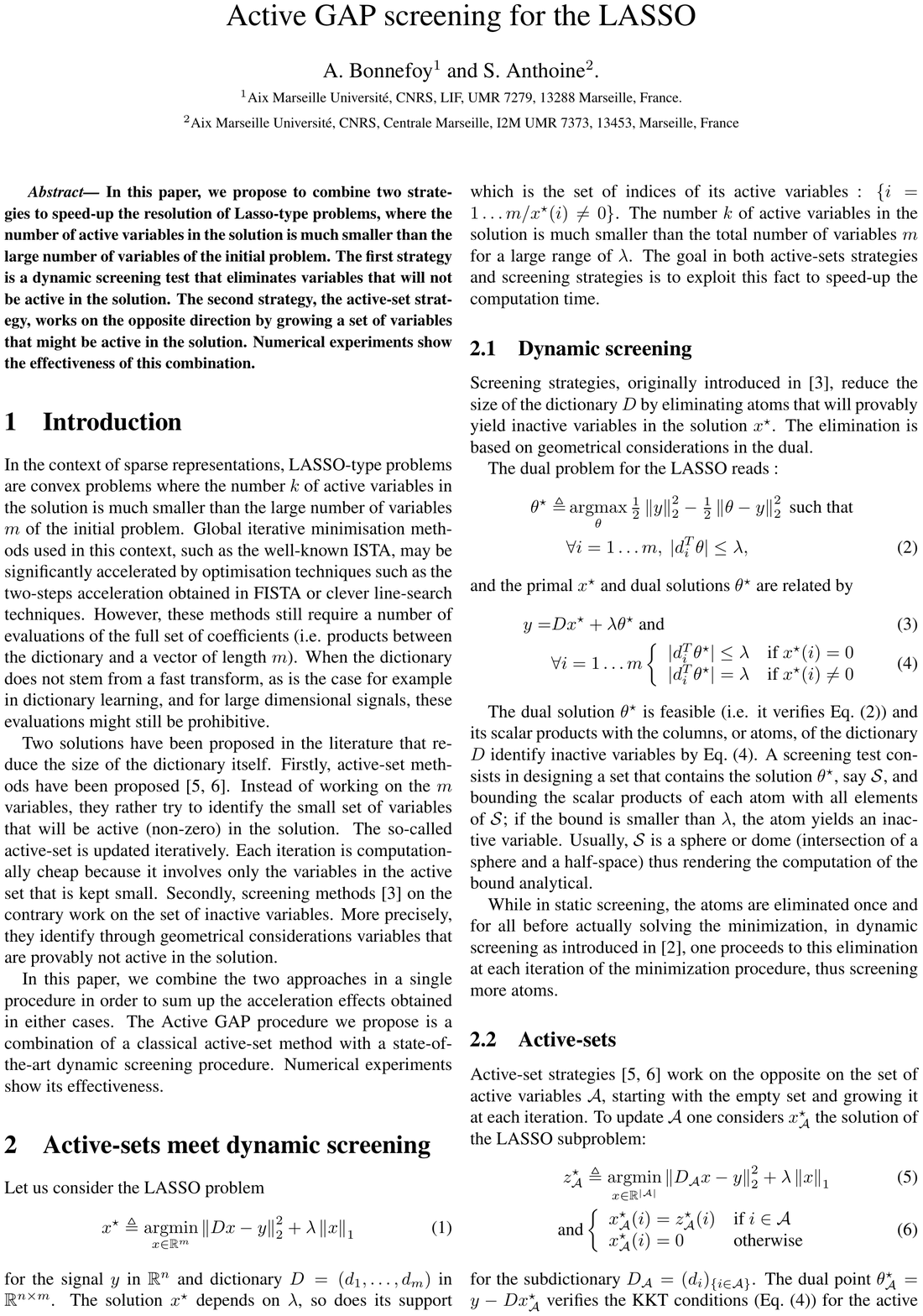}

\label{pdf:main_24}
\includepdf[offset=.65cm -2mm,pages=-,link=true,linkname=main_24,pagecommand={}]{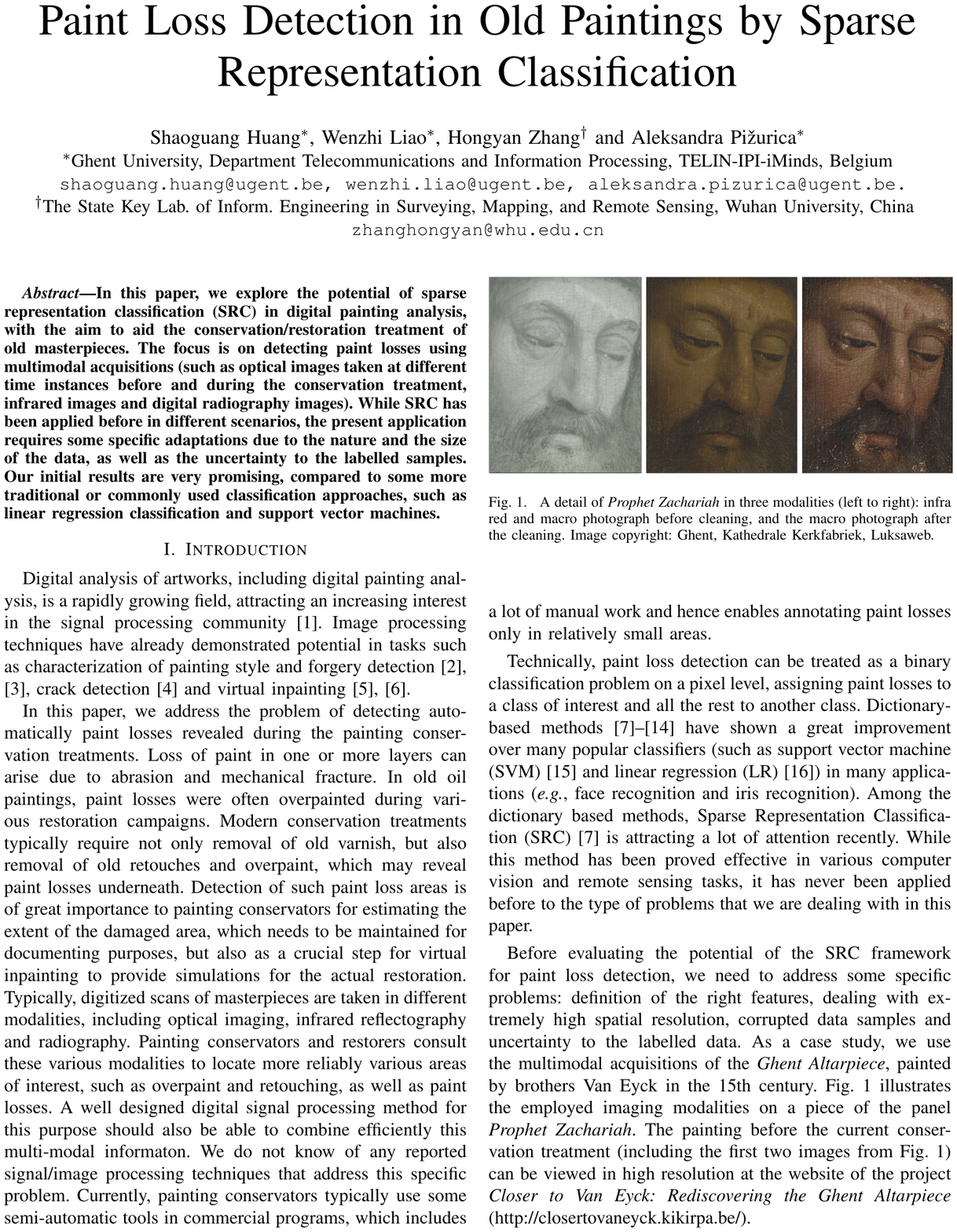}

\end{document}